\documentclass[smallcondensed]{svjour3_arXiv}     
\smartqed
\usepackage{graphicx,subfig}
\usepackage{amsmath,amssymb,amsfonts}
\usepackage{bm}
\usepackage{algorithm, algorithmic}
\usepackage{multirow,makecell}
\newcommand\gap{\hspace{0.1cm}}

\newcommand{\bS}{\mathbf{S}}
\newcommand{\N}{\mathcal{N}}
\newcommand{\T}{\mathcal{T}}
\newcommand{\Q}{\mathcal{Q}}

\newcommand{\un}{u^{(n)}}
\newcommand{\unn}{u^{(n+1)}}
\newcommand{\vn}{v^{(n)}}
\newcommand{\vnn}{v^{(n+1)}}
\newcommand{\Rw}{R_k^* w_k}

\newcommand{\intO}{\int_{\Omega}}
\newcommand{\sumk}{\sum_{k=1}^{N}}

\let\div\relax
\DeclareMathOperator{\div}{div}

\DeclareMathOperator{\dom}{dom}
\DeclareMathOperator*{\argmin}{\arg\min}

\numberwithin{equation}{section}

\counterwithin{theorem}{section}

\spnewtheorem{proposition}[theorem]{Proposition}{\bfseries}{\itshape}
\spnewtheorem{remark}[theorem]{Remark}{\itshape}{\rmfamily}

\usepackage{algorithm, algorithmic}

\begin{document}

\title{Accelerated Additive Schwarz Methods for Convex Optimization with Adaptive Restart
\thanks{This research was supported by Basic Science Research Program through the National Research Foundation of Korea~(NRF) funded by the Ministry of Education~(2019R1A6A1A10073887).}}
\titlerunning{Accelerated Additive Schwarz}
\author{Jongho Park}
\institute{Jongho Park \at
              Natural Science Research Institute, KAIST, Daejeon 34141, Korea \\
              \email{jongho.park@kaist.ac.kr}
}

\date{}

\maketitle

\begin{abstract}
Based on an observation that additive Schwarz methods for general convex optimization can be interpreted as gradient methods, we propose an acceleration scheme for additive Schwarz methods.
Adopting acceleration techniques developed for gradient methods such as momentum and adaptive restarting, the convergence rate of additive Schwarz methods is greatly improved.
The proposed acceleration scheme does not require any a priori information on the levels of smoothness and sharpness of a target energy functional, so that it can be applied to various convex optimization problems.
Numerical results for linear elliptic problems, nonlinear elliptic problems, nonsmooth problems, and nonsharp problems are provided to highlight the superiority and the broad applicability of the proposed scheme.
\end{abstract}

\section{Introduction}
\label{Sec:Introduction}
This paper is concerned with additive Schwarz methods for convex optimization problems.
Additive Schwarz methods are popular numerical solvers for large-scale linear elliptic problems, specialized for massively parallel computation; one may refer to~\cite{TW:2005,Xu:1992} for abstract theories of Schwarz methods for linear elliptic problems.
Meanwhile, there have been several successful applications of Schwarz methods to nonlinear problems; see, e.g.,~\cite{Badea:2006,BK:2012,Park:2020,TX:2002}.

An important observation on the Schwarz alternating method for linear elliptic problems is that the method can be viewed as a preconditioned Richardson method~\cite{TW:2005}.
Replacing Richardson iterations by conjugate gradient iterations with the same preconditioner, an improved algorithm that converges faster and does not require spectral information of linear operators is obtained.
This idea generalizes to general Schwarz methods for linear problems, so that most of modern Schwarz methods for linear problems are based on either the conjugate gradient method or the GMRES method.

An analogy of the Richardson method corresponding to general convex optimization is the gradient method.
Because of its simplicity and efficiency, there has been extensive research on the gradient method; see~\cite{CP:2016} for a survey on recent gradient methods and their applications.
In particular, starting from a celebrating work of Nesterov~\cite{Nesterov:1983},
acceleration of gradient methods has become an important topic in the field of mathematical optimization.
Momentum acceleration for smooth convex optimization was proposed in~\cite{Nesterov:1983}, and then improved in~\cite{DT:2014,KF:2016} to have better worst-case convergence rates.
In~\cite{BT:2009,CD:2015,Nesterov:2013}, it was generalized to nonsmooth convex optimization.
For strongly convex problems, it was shown in~\cite{CP:2016,Nesterov:2013,Nesterov:2018} that  further improvement of the convergence rate can be achieved by choosing the momentum value adaptively according to the level of strong convexity.
Alternatively, restarting techniques were adopted to deal with the strong convexity~\cite{Nesterov:2013,OC:2015,RD:2020}.

In the author's previous work~\cite{Park:2020}, it was proven that additive Schwarz methods for general convex optimization are interpreted as gradient methods.
In this perspective, abovementioned works on accelerated gradient methods may be considered for the sake of designing novel additive Schwarz methods that converge faster than existing ones.
Relevant existing works are~\cite{LPP:2019,LP:2019b}, in which accelerated domain decomposition methods for total variation minimization were obtained by FISTA acceleration~\cite{BT:2009}.
In this paper, we propose an acceleration scheme for additive Schwarz methods for convex optimization.
Integrating the plain additive Schwarz method with the Nesterov's momentum~\cite{Nesterov:1983} and the adaptive gradient restart scheme~\cite{OC:2015}, an accelerated method is obtained.
Differently from other acceleration schemes to deal with the sharpness~(see~\eqref{sharp} for the definition of the sharpness) of a target functional~\cite{CP:2016,Nesterov:2013,RD:2020}, the adaptive gradient restart scheme proposed in~\cite{OC:2015} does not require any information on the level of sharpness of the functional.
Therefore, the proposed scheme is applicable to a very broad range of convex optimization problems.
We present applications of the proposed scheme to additive Schwarz methods for nonlinear elliptic problems~\cite{TX:2002}, nonsmooth problems~\cite{Badea:2006,BTW:2003,Tai:2003,THX:2002}, and nonsharp problems~\cite{CTWY:2015,Park:2021}.
For all of those problems, we verify by numerical results that the proposed accelerated additive Schwarz methods outperform their unaccelerated counterparts.

The remainder of this paper is organized as follows.
In Section~\ref{Sec:ASM}, we briefly summarize key features of basic additive Schwarz methods for convex optimization.
We describe the proposed accelerated additive Schwarz method for convex optimization in Section~\ref{Sec:Accel}.
Applications of the proposed method to various convex optimization problems are presented in Section~\ref{Sec:Numerical}.
We conclude the paper with remarks in Section~\ref{Sec:Conclusion}.

\section{Additive Schwarz method}
\label{Sec:ASM}
In this section, we present a basic additive Schwarz method for the general convex optimization problem
\begin{equation}
\label{model}
\min_{u \in V} \left\{ E(u) := F(u) + G(u) \right\},
\end{equation}
where $V$ is a reflexive Banach space, $F$:~$V \rightarrow \mathbb{R}$ is a Frech\'{e}t differentiable convex functional, and $G$:~$V \rightarrow \overline{\mathbb{R}}$ is a proper, convex, lower semicontinuous functional that is possible nonsmooth.
We assume that $E$ is coercive so that there exists a solution $u^* \in V$ of~\eqref{model}.

Let $V_1$, $V_2$, \dots, and $V_N$ be reflexive Banach spaces.
In what follows, an index $k$ runs from 1 to $N$.
We assume that there exist a bounded linear operator $R_k^*$:~$V_k \rightarrow V$ such that
\begin{equation}
\label{space_decomp}
V = \sumk R_k^* V_k
\end{equation}
and its adjoint $R_k$:~$V \rightarrow V_k$ is surjective.
Under the space decomposition~\eqref{space_decomp}, an additive Schwarz method for~\eqref{model} with exact local solvers is given as Algorithm~\ref{Alg:ASM}.
For more general setting that allows inexact local solvers, see~\cite{Park:2020}.

\begin{algorithm}[]
\caption{Additive Schwarz method for~\eqref{model}}
\begin{algorithmic}[]
\label{Alg:ASM}
\STATE Let $u^{(0)} \in \dom G$ and $\tau > 0$.
\FOR{$n=0,1,2,\dots$}
\item \begin{equation*}
\begin{split}
w_k^{(n+1)} &\in \argmin_{w_k \in V_k} E(\un + R_k^* w_k ), \quad 1 \leq k \leq N \\
\unn &= \un + \tau \sumk \Rw^{(n+1)}
\end{split}
\end{equation*}
\ENDFOR
\end{algorithmic}
\end{algorithm} 
 
As a special case of Algorithm~\ref{Alg:ASM}, we consider the case of linear elliptic problems.
Suppose temporarily that $V$, $V_k$ are Hilbert spaces and that the energy functional $E$ in~\eqref{model} is given by
\begin{equation}
\label{linear_model}
F(u) = \frac{1}{2}\left<Au, u \right> - \left< f, u \right>, \quad G(u) = 0,
\end{equation}
where $A$:~$V \rightarrow V$ is a continuous and symmetric positive definite linear operator and $f \in V$.
We define the local stiffness operator $A_k$:~$V_k \rightarrow V_k$ by
\begin{equation*}
\left< A_k u_k, v_k \right> = \left<A R_k^* u_k, R_k^* v_k \right>, \quad u_k, v_k \in V_k.
\end{equation*}
Then it is straightforward to show that~(see, e.g.,~\cite[section~4.1]{Park:2020}) Algorithm~\ref{Alg:ASM} for~\eqref{linear_model} can be rewritten as
\begin{subequations}
\label{Richardson}
\begin{equation}
\unn = \un - \tau M^{-1} (A\un - f), \quad n \geq 0,
\end{equation}
or equivalently,
\begin{equation}
\resizebox{\textwidth}{!}{$ \displaystyle
\unn = \argmin_{u \in V} \left\{ F(\un) + \langle F'(\un), u-\un \rangle + \frac{1}{2\tau}  \langle M(u - \un), u-\un \rangle \right\}, \quad n \geq 0,
$}
\end{equation}
\end{subequations}
where $M$:~$V \rightarrow V$ is the additive Schwarz preconditioner given by
\begin{equation}
\label{preconditioner}
M = \left( \sumk R_k^* A_k^{-1} R_k \right)^{-1}.
\end{equation}
Equation~\eqref{Richardson} implies that Algorithm~\ref{Alg:ASM} for~\eqref{linear_model} is the Richardson method for the preconditioned system
\begin{equation}
\label{preconditioned}
M^{-1}Au = M^{-1}f.
\end{equation}
Therefore, by applying the conjugate gradient method to~\eqref{preconditioned} instead of the Richardson method, we can obtain a more improved algorithm than Algorithm~\ref{Alg:ASM}.
For a theoretical comparison of the conjugate gradient method with the Richardson method, one may refer to~\cite[Appendix~C]{TW:2005}.

In~\cite[Lemma~4.5]{Park:2020}, it was observed that~\eqref{Richardson} can be generalized to additive Schwarz methods for the general convex optimization~\eqref{model}.
A rigorous statement is presented in the following proposition.

\begin{proposition}[generalized additive Schwarz lemma]
\label{Prop:ASM}
Let $\{ \un \}$ be the sequence generated by Algorithm~\ref{Alg:ASM}.
Then it satisfies
\begin{equation}
\label{GD}
\unn \in \argmin_{u \in V} \left\{ F(\un) + \langle F'(\un), u - \un \rangle + M_{\tau} (u, u^{(n)}) \right\}, \quad n \geq 0,
\end{equation}
where the functional $M_{\tau}$:~$V \times V \rightarrow \overline{\mathbb{R}}$ is given by
\begin{equation*}
\begin{split}
M_{\tau} (u, v) &= \tau \inf \left\{ \sumk ( \left( D_F(v + \Rw , v) + G ( v + \Rw) \right)  : u-v = \tau \sumk \Rw, \gap w_k \in V_k \right\} \\
 &\quad+ \left(1- \tau N \right) G( v), \quad u,v \in V, 
\end{split}
\end{equation*}
and $D_F$:~$V \times V \rightarrow \mathbb{R}$ is the Bregman distance of $F$ defined by
\begin{equation*}
D_F (u,v) = F(u) - F(v) - \left< F'(v), u-v \right>, \quad u,v \in V.
\end{equation*}
\end{proposition}

It is clear that~\eqref{GD} reduces to~\eqref{Richardson} in the case of~\eqref{linear_model}.
Proposition~\ref{Prop:ASM} means that Algorithm~\ref{Alg:ASM} is an instance of nonlinear gradient methods for~\eqref{model}; see~\cite{Teboulle:2018} for a recent survey on nonlinear gradient methods of the form~\eqref{GD}.
In~\cite{Park:2020}, an abstract convergence theory of additive Schwarz methods that generalizes~\cite[Chapter~2]{TW:2005} was developed using Proposition~\ref{Prop:ASM} and the convergence theory of nonlinear gradient methods.

\section{Acceleration schemes}
\label{Sec:Accel}
First, we review existing acceleration schemes for gradient methods for~\eqref{model}.
We recall that the energy functional $E$ is said to be sharp if there exists a constant $p > 1$ such that for any bounded and convex subset $K$ of $V$ satisfying $u^* \in K$, we have
\begin{equation}
\label{sharp}
\frac{\mu_K}{p} \| u - u^* \|^p \leq E(u) - E(u^*), \quad u\in K,
\end{equation}
for some $\mu_K > 0$.

As a fundamental example of gradient methods, we consider the forward-backward splitting method~\cite{BT:2009,CW:2005}, also known as the composite gradient method~\cite{Nesterov:2013}.
In~\eqref{model}, assume that $F'$ is Lipschitz continuous with modulus $L$, i.e., it satisfies
\begin{equation*}
F(u) \leq F(v) + \left< F'(v), u-v \right> + \frac{L}{2} \| u - v\|^2, \quad u,v \in V.
\end{equation*}
The forward-backward splitting method for~\eqref{model} is presented in Algorithm~\ref{Alg:FB}.

\begin{algorithm}[]
\caption{Forward-backward splitting for~\eqref{model}}
\begin{algorithmic}[]
\label{Alg:FB}
\STATE Let $u^{(0)} \in \dom G$ and $\tau \in (0, 1/L]$.
\FOR{$n=0,1,2,\dots$}
\item \begin{equation*}
\unn \in \argmin_{u \in V} \left\{ F(\un) + \langle F'(\un), u - \un \rangle + \frac{1}{2\tau} \| u - \un \|^2 + G(u) \right\}
\end{equation*} 
\ENDFOR
\end{algorithmic}
\end{algorithm}

It is well-known that the worst-case energy error of Algorithm~\ref{Alg:FB} decays with the rate $O(1/n)$~\cite{BT:2009,Nesterov:2013}.
If the energy functional $E$ is sharp, then an improved error bound can be obtained~\cite{Park:2020,RD:2020}.
In particular, under the assumptions that $E$ is strongly convex, i.e., when it satisfies~\eqref{sharp} with $p = 2$, Algorithm~\ref{Alg:FB} converges linearly. 

Algorithm~\ref{Alg:FB} can be accelerated by adding momentum.
At each step of the algorithm, we set
\begin{equation*}
\vnn = \unn + \beta_n (\unn - \un)
\end{equation*}
for some suitably chosen $\beta_n > 0$, and then apply the forward-backward splitting to $\vn$ instead of $\un$ in the next step.
Such an acceleration scheme was first proposed by Nesterov~\cite{Nesterov:1983} for smooth convex optimization~($G=0$ in~\eqref{model}), and then generalized to the nonsmooth case in~\cite{BT:2009}.
Among several variants of the momentum technique~\cite{BT:2009,CD:2015,Nesterov:2013} for~\eqref{model}, we present FISTA~\cite{BT:2009} in Algorithm~\ref{Alg:FISTA}.

\begin{algorithm}[]
\caption{FISTA for~\eqref{model}}
\begin{algorithmic}[]
\label{Alg:FISTA}
\STATE Let $u^{(0)} = v^{(0)} \in \dom G$, $\tau \in (0, 1/L]$, and $t_0 = 1$.
\FOR{$n=0,1,2,\dots$}
\item \begin{equation*}
\unn \in \argmin_{u \in V} \left\{ F(\vn) + \langle F'(\vn), u - \vn \rangle + \frac{1}{2\tau} \| u - \vn \|^2 + G(u) \right\}
\end{equation*} 
\item \begin{equation*}
t_{n+1} = \frac{1 + \sqrt{1 + 4t_n^2}}{2}, \quad
\beta_n = \frac{t_n - 1}{t_{n+1}}
\end{equation*}
\item \begin{equation*}
\vnn = \unn + \beta_n (\unn - \un)
\end{equation*}
\ENDFOR
\end{algorithmic}
\end{algorithm}

It was shown in~\cite[Theorem~4.4]{BT:2009} that Algorithm~\ref{Alg:FISTA} enjoys the $O(1/n^2)$ convergence rate, which is faster than Algorithm~\ref{Alg:FB}.
This rate is optimal for smooth convex optimization in the sense that there exists a smooth convex functional such that any first-order method for minimizing the functional must satisfy an $O(1/n^2)$ lower bound of the energy error; see, e.g.,~\cite[Theorem~4.3]{CP:2016}.
However, when the energy functional $E$ is sharp, Algorithm~\ref{Alg:FISTA} is not enough to guarantee the optimal convergence rate; the momentum parameter $\beta_n$ must be chosen according to the sharpness information of $E$.
In~\cite{CP:2016,Nesterov:2013}, momentum techniques suitable for the strongly convex case were considered.
Alternatively, the optimal rate can be achieved by restarting Algorithm~\ref{Alg:FISTA} appropriately; we reset the momentum parameters as $t_{n+1} = 1$ and $\beta_n = 0$ whenever the iterates of the algorithm meet some criterion.
A restarting technique for the strongly convex objective functional was considered in~\cite{Nesterov:2013}, and then it was generalized to the general sharp case in~\cite{RD:2020}.
All of the abovementioned approaches to deal with the sharp case share a common drawback that they require explicit values for the sharpness information $p$ and $\mu_K$ in~\eqref{sharp}.
Since a priori sharpness information of the energy functional is not available in general, such a drawback is crucial in practice.
In~\cite{OC:2015}, adaptive restarting techniques were proposed which are heuristic but very effective acceleration schemes.
Although they do not require any information on the sharpness of the energy functional, it was numerically verified that their performances are as good as the abovementioned methods.
Algorithm~\ref{Alg:FISTA_ada} presents the gradient adaptive restart scheme proposed in~\cite{OC:2015}, applied to Algorithm~\ref{Alg:FISTA}.
In the criterion for restart in Algorithm~\ref{Alg:FISTA_ada}, $E'(\vn)$ denotes the composite gradient~\cite{Nesterov:2013} of $E$ at $\vn$, a notion that generalizes the usual gradient for composite objective functionals of the form~\eqref{model}.

\begin{algorithm}[]
\caption{FISTA with adaptive restart for~\eqref{model}}
\begin{algorithmic}[]
\label{Alg:FISTA_ada}
\STATE Let $u^{(0)} = v^{(0)} \in \dom G$, $\tau \in (0, 1/L]$, and $t_0 = 1$.
\FOR{$n=0,1,2,\dots$}
\item \begin{equation*}
\unn = \argmin_{u \in V} \left\{ F(\vn) + \langle F'(\vn), u - \vn \rangle + \frac{1}{2\tau} \| u - \vn \|^2 + G(u) \right\}
\end{equation*}
\IF {$\left< E'( \vn ) , \unn - \un \right> > 0$}
\STATE \begin{equation*}
t_{n+1} = 1, \quad
\beta_n = 0
\end{equation*}
\ELSE
\STATE \begin{equation*}
t_{n+1} = \frac{1 + \sqrt{1 + 4t_n^2}}{2}, \quad
\beta_n = \frac{t_n - 1}{t_{n+1}}
\end{equation*}
\ENDIF
\STATE \begin{equation*}
\vnn = \unn + \beta_n (\unn - \un)
 \end{equation*}
\ENDFOR
\end{algorithmic}
\end{algorithm}

Now, we are ready to propose an accelerated additive Schwarz method for~\eqref{model}.
Combining the additive Schwarz method presented in Algorithm~\ref{Alg:ASM} with the idea of gradient adaptive restart, we propose an accelerated version of Algorithm~\ref{Alg:ASM}; see Algorithm~\ref{Alg:proposed}.

\begin{algorithm}[]
\caption{Accelerated additive Schwarz method for~\eqref{model}}
\begin{algorithmic}[]
\label{Alg:proposed}
\STATE Let $u^{(0)} = v^{(0)} \in \dom G$, $\tau > 0$, and $t_0 = 1$.
\FOR{$n=0,1,2,\dots$}
\item \begin{equation*}
\begin{split}
w_k^{(n+1)} &\in \argmin_{w_k \in V_k} E(\vn + R_k^* w_k ), \quad 1 \leq k \leq N \\
\unn &= \vn + \tau \sumk \Rw^{(n+1)}
\end{split}
\end{equation*}
\IF {$\left< \vn - \unn, \unn - \un \right> > 0$}
\STATE \begin{equation*}
t_{n+1} = 1, \quad
\beta_n = 0
\end{equation*} 
\ELSE
\STATE \begin{equation*}
t_{n+1} = \frac{1 + \sqrt{1 + 4t_n^2}}{2}, \quad
\beta_n = \frac{t_n - 1}{t_{n+1}}
\end{equation*} 
\ENDIF
\STATE \begin{equation*}
\vnn = \unn + \beta_n (\unn - \un)
 \end{equation*} 
\ENDFOR
\end{algorithmic}
\end{algorithm}

In Algorithm~\ref{Alg:proposed}, we see that
\begin{equation*}
\vn - \unn = - \tau \sumk R_k^* w_k^{(n+1)}.
\end{equation*}
In view of Proposition~\ref{Prop:ASM}, one may regard the right-hand side of the above equation as a ``generalized'' gradient of the energy functional $E$ at $\vn$ with respect to the non-Euclidean distance function $M_{\tau}$.
In this sense, we replace the composite gradient $E'(\vn)$ in the restart criterion of Algorithm~\ref{Alg:FISTA_ada} by $\vn - \unn$ in the proposed method.
Since Proposition~\ref{Prop:ASM} means that the plain additive Schwarz method presented in Algorithm~\ref{Alg:ASM} is the gradient method for~\eqref{model} with respect to $M_{\tau}$, it is expected that the restarting step in Algorithm~\ref{Alg:proposed} can improve the convergence rate of the additive Schwarz method by the same principle as Algorithm~\ref{Alg:FISTA_ada}.
More precisely, the restart criterion
\begin{equation}
\label{restart_criterion}
\left< \vn - \unn, \unn - \un \right> > 0
\end{equation}
in Algorithm~\ref{Alg:proposed} means that the update direction $\unn - \un$ is on the same side of the generalized gradient direction $\vn - \unn$.
Since the energy decreases toward the minus gradient direction in general, meeting the restart criterion~\eqref{restart_criterion} implies that the overrelaxed variable $\vn = \un + \beta_{n-1}(\un - u^{(n-1)})$ was not properly chosen.
Hence, it is natural to consider resetting the overrelaxation parameter $\beta_n$ as 0 whenever~\eqref{restart_criterion} is satisfied.

\begin{remark}
\label{Rem:function}
In~\cite{OC:2015}, the function adaptive restart scheme that restarts FISTA whenever $E(\unn) > E(\un)$ was proposed as well as the gradient adaptive restart scheme.
As an alternative of Algorithm~\ref{Alg:proposed}, one may adopt the function adaptive restart scheme for additive Schwarz method in order to accelerate the convergence.
However, the function adaptive restart scheme has a disadvantage that additional computational cost for $E(\unn)$ is need in each iteration.
On the contrary, Algorithm~\ref{Alg:proposed} does not require additional major computational cost since $\un$, $\unn$, and $\vn$ are computed prior to checking the restart criterion.
In this perspective, we do not deal with the function adaptive restart scheme in this paper.
Similar discussions were made in~\cite{OC:2015}.
\end{remark}

\begin{figure}[]
\centering
\subfloat[][]{ \includegraphics[height=0.31\linewidth]{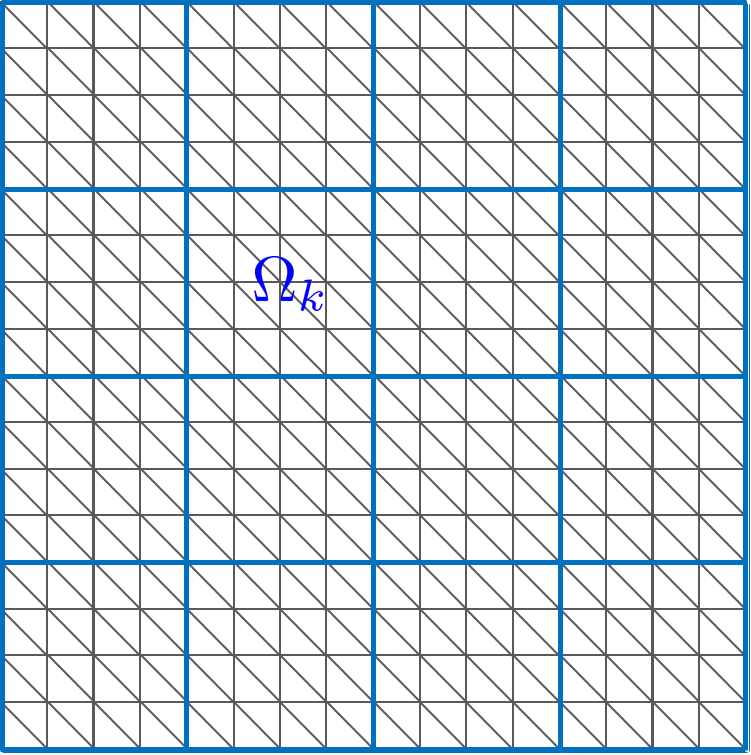} }
\gap
\subfloat[][]{ \includegraphics[height=0.31\linewidth]{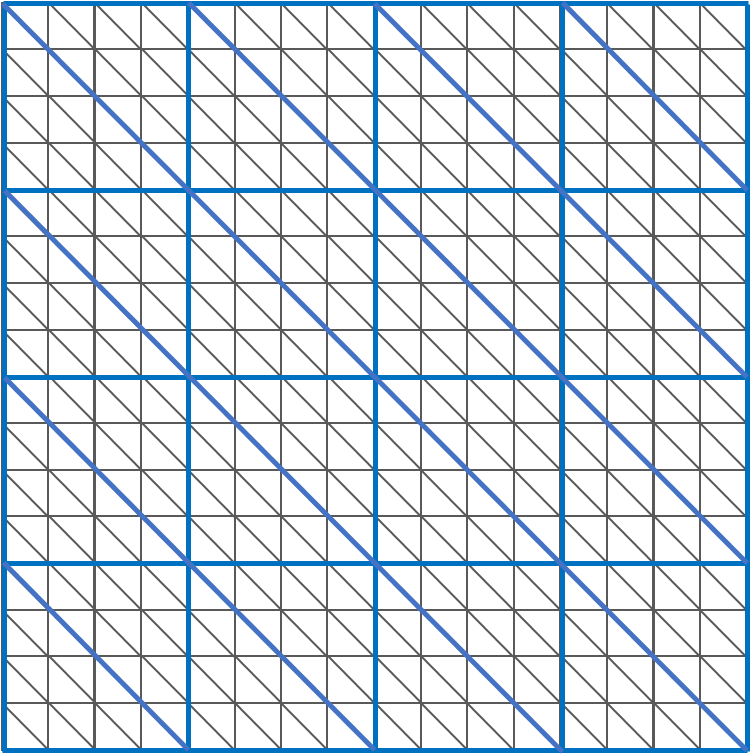} }
\gap
\subfloat[][]{ \includegraphics[height=0.31\linewidth]{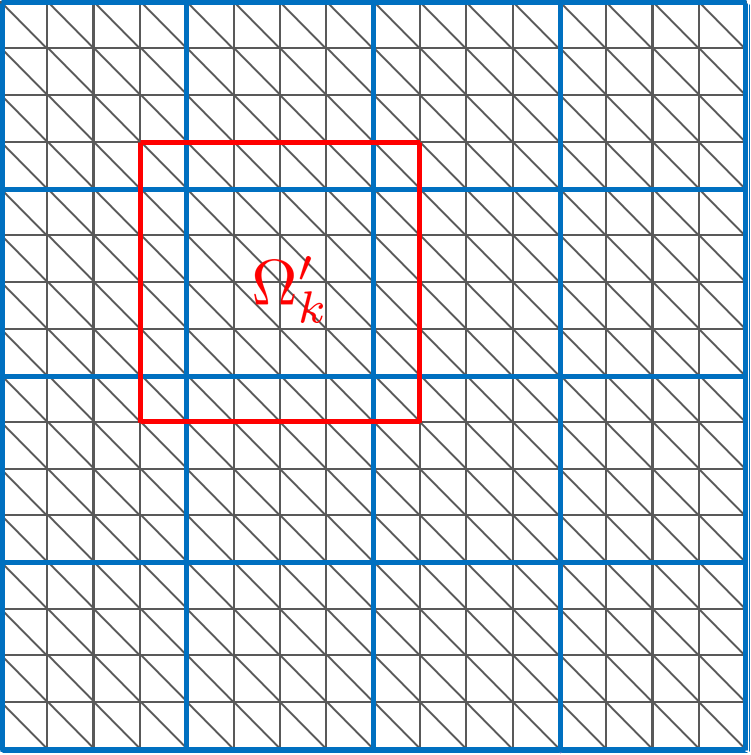} }
\caption{Domain decomposition setting when $h = 1/2^4$, $H = 1/2^2$, and $\delta = h$: \textbf{(a)}~nonoverlapping domain decomposition $\{ \Omega_k \}_{k=1}^{\N}$, \textbf{(b)}~coarse triangulation $\T_H$, \textbf{(c)}~overlapping domain decomposition $\{ \Omega_k' \}_{k=1}^{\N}$.}
\label{Fig:DD}
\end{figure}

The major part of each iteration of Algorithm~\ref{Alg:proposed} is to solve local minimization problems on $V_k$; the computation cost for momentum parameters $t_n$ and $\beta_n$ is clearly marginal.
Therefore, the main computational effort of Algorithm~\ref{Alg:proposed} is the same as the one of Algorithm~\ref{Alg:ASM}.
In addition, inheriting the advantage of Algorithm~\ref{Alg:FISTA_ada}, the proposed method does not require any prior information on the levels of smoothness and sharpness of the energy functional $E$.
Choosing the step size $\tau$ of Algorithm~\ref{Alg:proposed} usually depends on a domain decomposition setting but not on the energy functional.
For example, in the usual one- and two-level overlapping domain decomposition settings for a two-dimensional domain~(see Figure~\ref{Fig:DD}(c)), one may set $\tau = 1/4$ and $1/5$, respectively, since the subdomains can be colored with 4 colors; see~\cite[Section~5.1]{Park:2020} for details.

\section{Numerical experiments}
\label{Sec:Numerical}
In this section, we present applications of Algorithm~\ref{Alg:proposed} to various nonlinear problems   appearing in science and engineering that can be represented in the form~\eqref{model}.
In particular, we conducted numerical experiments on nonlinear elliptic problem, obstacle problem, and dual total variation minimization.
For all the problems, we claim that the proposed method has a superior convergence property compared to the unaccelerated one.
All computations presented in this section were performed on a computer cluster equipped with Intel Xeon SP-6148 CPUs~(2.4GHz, 20C) and the operating system CentOS~7.4 64bit. 

\subsection{Nonlinear elliptic problem}
\label{Subsec:sLap}
First, we present an application of the proposed method to the following model $s$-Laplace equation:
\begin{equation} \begin{split}
\label{sLap_PDE} 
- \div \left( |\nabla u|^{s-2} \nabla u \right) &= f \quad\textrm{ in } \Omega , \\
u &= 0 \quad\textrm{ on } \partial \Omega ,
\end{split} \end{equation}
where $s > 1$ and $f \in L^{\frac{s}{s-1}}(\Omega)$.
We note that Schwarz methods for the problem~\eqref{sLap_PDE} were considered in~\cite{TX:2002}.
It is well-known that~(see~\cite{Ciarlet:2002} for instance) a unique solution of~\eqref{sLap_PDE} solves the following minimization problem:
\begin{equation}
\label{sLap_opt}
\min_{u \in W_0^{1,s}(\Omega)} \left\{ \frac{1}{s} \intO |\nabla u|^s \,dx - \intO fu \,dx \right\}.
\end{equation}

In the following, we set $\Omega = [0,1]^2 \subset \mathbb{R}^2$.
We decompose the domain $\Omega$ into $\N = N \times N$ square subdomains $\{ \Omega_k \}_{k=1}^{\N}$in which each subdomain has the sidelength $H=1/N$.
Each subdomain $\Omega_k$, $1\leq k \leq \N$, is partitioned into $2 \times H/h \times H/h$ uniform triangles to form a global triangulation $\T_h$ of $\Omega$.
Similarly, we partition each $\Omega_k$ into two uniform triangles and let $\T_H$ be a coarse triangulation of $\Omega$ consisting of such triangles.
Overlapping subdomains $\{ \Omega_k' \}_{k=1}^{\N}$ are constructed in a way that $\Omega_k'$ is a union of $\Omega_k$ and its surrounding layers of fine elements in $\T_h$ with the width $\delta$ such that $0< \delta < H/2$.
Figure~\ref{Fig:DD} illustrates the domain decomposition of $\Omega$ explained above.

Let $S_h (\Omega)$ and $S_H (\Omega)$ be the $\mathcal{P}_1$-Lagrangian finite element spaces on $\T_h$ and $\T_H$ with the homogeneous essential boundary condition, respectively.
A conforming approximation of~\eqref{sLap_opt} using $S_h (\Omega) \subset W_0^{1,s} (\Omega)$ is written as
\begin{equation}
\label{sLap_FEM}
\min_{u \in S_h (\Omega)} \left\{ \frac{1}{s} \intO |\nabla u|^s \,dx - \intO fu \,dx \right\}.
\end{equation}
The discretized problem~\eqref{sLap_FEM} can be represented in the form of~\eqref{model} with
\begin{equation*}
V = S_h (\Omega), \quad F(u) = \frac{1}{s} \intO |\nabla u|^s \,dx - \intO fu \,dx, \quad G(u) = 0.
\end{equation*}
One can show that $F$ satisfies the sharpness condition~\eqref{sharp} with $p=s$ and $p=2$ when $s>2$ and $1<s<2$, respectively~\cite{Ciarlet:2002}.

If we set
\begin{equation*}
V_k = S_h (\Omega_k'), \quad 1 \leq k \leq \N,
\end{equation*}
and take $R_k^*$:~$V_k \rightarrow V$ as the natural extension operator, then it clearly satisfies the space decomposition assumption~\eqref{space_decomp}, where $S_h (\Omega_k')$ is the $\mathcal{P}_1$-Lagrangian finite element space on the $\T_h$-elements in $\Omega_k'$ with the homogeneous essential boundary condition.
For the two-level setting, we set
\begin{equation*}
V_0 = S_H (\Omega)
\end{equation*}
so that
\begin{equation}
\label{space_decomp_2L}
V = R_0^* V_0 + \sumk R_k^* V_k,
\end{equation}
where $R_0^*$:~$V_0 \rightarrow V$ is the natural interpolation operator.
Under the space decompositions~\eqref{space_decomp} and~\eqref{space_decomp_2L}, convergence of the unaccelerated additive Schwarz method for~\eqref{sLap_FEM} is guaranteed by the following proposition~\cite[Theorem~6.1]{Park:2020}.

\begin{proposition}
\label{Prop:sLap}
In Algorithm~\ref{Alg:ASM} for~\eqref{sLap_FEM}, if $E(u^{(0)}) - E(u^*)$ is small enough, then there exists a positive constant $C$ independent of $h$, $H$, and $\delta$ such that
\begin{equation*}
E(u^{(n)}) - E(u^*) \leq C \frac{1 + 1/\delta^q}{(n+1)^{\frac{p(q-1)}{p-q}}}
\end{equation*}
for the one-level domain decomposition~\eqref{space_decomp} with $\tau \leq 1/4$ and
\begin{equation*}
E(u^{(n)}) - E(u^*) \leq C \frac{1 + (H/\delta)^q}{(n+1)^{\frac{p(q-1)}{p-q}}}
\end{equation*}
for the two-level domain decomposition~\eqref{space_decomp_2L} with $\tau \leq 1/5$, where
\begin{align*}
p=s, \gap q=2 \quad &\textrm{ if } s>2, \\
p=2, \gap q=s \quad &\textrm{ if } 1<s<2.
\end{align*}
\end{proposition}

\begin{remark}
Proposition~\ref{Prop:sLap} implies that the unaccelerated additive Schwarz method for~\eqref{sLap_FEM} satisfies the $O(1/n^{\frac{p(q-1)}{p-q}})$ convergence rate.
We compare this estimate with existing ones~\cite{Badea:2006,BK:2012,TX:2002} to show that our estimate is sharper than the existing results; see also~\cite{Park:2020}.
The first rigorous analysis on the convergence rate of the unaccelerated additive Schwarz method for~\eqref{sLap_FEM} was presented in~\cite{TX:2002}; the $O(1/n^{\frac{q(q-1)}{(p-q)(p+q-1)}}$ convergence rate was proven.
More recently, the $O(1/n^{\frac{q-1}{p-q}})$ convergence rate was analyzed in~\cite{Badea:2006,BK:2012}.
Since
\begin{equation*}
\frac{q(q-1)}{(p-q)(p+q-1)} < \frac{q-1}{p-q} < \frac{p(q-1)}{p-q}
\end{equation*}
for $1<q<p$, one can conclude that Proposition~\ref{Prop:sLap} provides a sharper estimate than the existing results~\cite{Badea:2006,BK:2012,TX:2002}.
\end{remark}

Now, we present numerical results of the proposed method applied to~\eqref{sLap_FEM} with $s= 4$ and $f = 1$.
For all experiments, we set the initial guess as $u^{(0)} = 0$.
Local problems on $V_k$, $1\leq k \leq \N$, and coarse problems on $V_0$ were solved by Algorithm~\ref{Alg:FISTA_ada} equipped with the backtracking strategy proposed in~\cite{BT:2009} using the stop criteria
\begin{equation}
\label{stop_local}
h^2 \| w_k^{(n+1)} - w_k^{(n)} \|_{\ell^2}^2 < 10^{-20}
\end{equation}
and
\begin{equation}
\label{stop_coarse}
H^2 \| w_0^{(n+1)} - w_0^{(n)} \|_{\ell^2}^2 < 10^{-20}
\end{equation}
respectively, where $\| \cdot \|_{\ell^2}$ denotes the $\ell^2$-norm of degrees of freedom.
The step size $\tau$ in Algorithms~\ref{Alg:ASM} and~\ref{Alg:proposed} was chosen as $\tau=1/4$ for the one-level decomposition~\eqref{space_decomp} and $\tau=1/5$ for the two-level decomposition~\eqref{space_decomp_2L}.
A reference solution $u^* \in V$ was computed by $10^5$ iterations of Algorithm~\ref{Alg:FISTA_ada} for the full-dimension problem~\eqref{sLap_FEM}. 
 
\begin{figure}[]
\centering
\subfloat[][One-level decomposition~\eqref{space_decomp}]{ \includegraphics[width=0.47\linewidth]{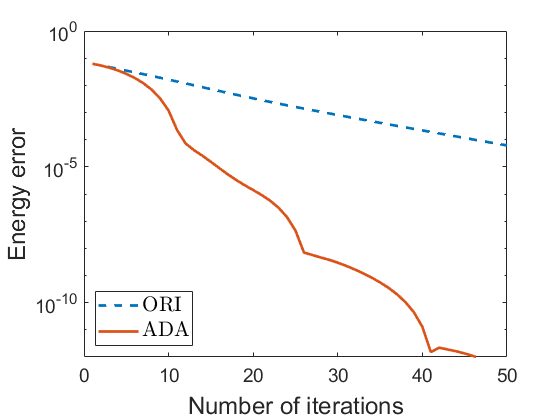} }
\gap
\subfloat[][Two-level decomposition~\eqref{space_decomp_2L}]{ \includegraphics[width=0.47\linewidth]{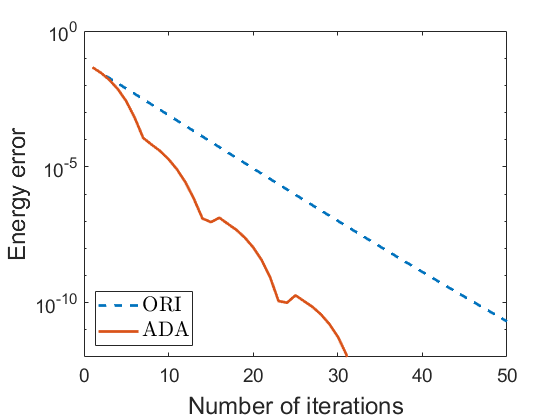} }
\caption{Decay of the energy error $E(u^{(n)}) - E(u^*)$ in additive Schwarz methods for the $s$-Laplacian problem~\eqref{sLap_FEM}~($h = 1/2^6$, $H=1/2^3$, $\delta=4h$).
ORI and ADA denote the unaccelerated~(Algorithm~\ref{Alg:ASM}) and accelerated~(Algorithm~\ref{Alg:proposed}) methods, respectively.}
\label{Fig:sLap}
\end{figure}

In order to highlight the efficiency of the proposed method for~\eqref{sLap_FEM}, we compare the energy decay of the unaccelerated and accelerated methods.
We note that the unaccelerated method, Algorithm~\ref{Alg:ASM} for~\eqref{sLap_FEM}, is identical to~\cite[Algorithm~2.1]{TX:2002}.
Figure~\ref{Fig:sLap} plots the energy error $E(u^{(n)}) - E(u^*)$ of Algorithms~\ref{Alg:ASM} and~\ref{Alg:proposed} when $h=1/2^6$, $H=1/2^3$, and $\delta = 4h$.
For both of the cases one-level and two-level domain decomposition settings, the proposed method shows faster convergence to the energy minimum compared to the unaccelerated method.
Since the main computational costs of Algorithms~\ref{Alg:ASM} and~\ref{Alg:proposed} are the same, we can say that the proposed method is superior to the conventional method in the sense of both convergence rate and computational cost.

Proposition~\ref{Prop:sLap} implies that Algorithm~\ref{Alg:ASM} is scalable in the sense that its convergence rate depends only on the size of local problems $H/h$ whenever $\delta/h$ is fixed.
That is, when each subdomain is assigned to a single processor, Algorithm~\ref{Alg:ASM} can solve a problem of the larger size with the same amount of time if more parallel processors can be utilized simultaneously.
Since the proposed method showed superior convergence results compared to Algorithm~\ref{Alg:ASM} in the above numerical experiments, one can readily expect that it is also scalable.
In the following, we verify the scalability of Algorithm~\ref{Alg:proposed} by numerical experiments.

\begin{figure}[]
\centering
\subfloat[][$H/h = 2^3$, $\delta/h = 2$]{ \includegraphics[width=0.47\linewidth]{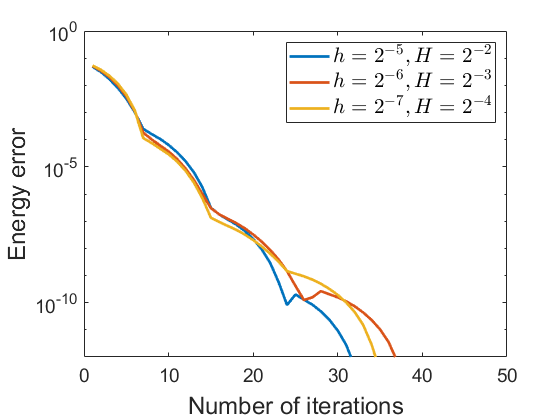} }
\gap
\subfloat[][$H/h = 2^4$, $\delta/h = 2$]{ \includegraphics[width=0.47\linewidth]{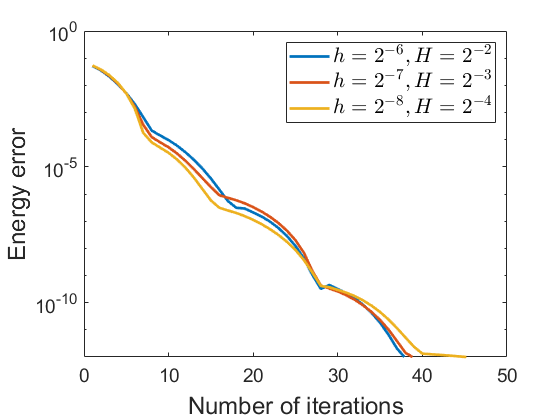} }
\caption{Decay of the energy error $E(\un) - E(u^*)$ in Algorithm~\ref{Alg:proposed} for the $s$-Laplacian problem~\eqref{sLap_FEM} when $H/h$ and $\delta/h$ are fixed.}
\label{Fig:sLap_sca}
\end{figure}

\begin{table}
\centering
\begin{tabular}{c c c c} \hline
$H/h$ & $H$ & $h$ & \#iter \\
\Xhline{1pt}
\multirow{3}{*}{$2^2$} & $1/2^2$ & $1/2^4$ & 20 \\
& $1/2^3$ & $1/2^5$ & 21 \\
& $1/2^4$ & $1/2^6$ & 20 \\
\hline
\multirow{3}{*}{$2^3$} & $1/2^2$ & $1/2^5$ & 21 \\
& $1/2^3$ & $1/2^6$ & 22 \\
& $1/2^4$ & $1/2^7$ & 22 \\
\hline
\multirow{3}{*}{$2^4$} & $1/2^2$ & $1/2^6$ & 26 \\
& $1/2^3$ & $1/2^7$ & 26 \\
& $1/2^4$ & $1/2^8$ & 25 \\
\hline
\end{tabular}
\caption{Number of iterations of Algorithm~\ref{Alg:proposed} for the $s$-Laplacian problem~\eqref{sLap_FEM} to meet the condition~\eqref{stop_sca}.}
\label{Table:sLap_sca}
\end{table}

Figure~\ref{Fig:sLap_sca} displays the energy decay of Algorithm~\ref{Alg:proposed} for~\eqref{sLap_FEM} when $H/h$ and $\delta/h$ are fixed.
One can observe that the convergence rate of Algorithm~\ref{Alg:proposed} remains almost the same when both $h$ and $H$ decrease keeping their ratio $H/h$ constant.
In addition, as shown in Table~\ref{Table:sLap_sca}, the numbers of iterations to satisfy the stop condition
\begin{equation}
\label{stop_sca}
E(\un) - E(u^*) < 10^{-8}
\end{equation}
for various $H$ and $h$ are uniformly bounded for a fixed value of $H/h$.
It shows the numerical scalability of Algorithm~\ref{Alg:proposed}.

\subsection{Obstacle problem}
Next, we apply the proposed method to the following model variational inequality: find $u \in K$ such that
\begin{equation}
\label{obs_VI}
\intO \nabla u \cdot \nabla (v-u) \,dx \geq 0 \quad \forall v \in K,
\end{equation}
where $K$ is a nonempty convex subset of $H_0^1 (\Omega)$ defined in terms of obstacle functions $g_L, g_U \in L^{\infty}(\Omega)$:
\begin{equation*}
K = \left\{ u \in H_0^1 (\Omega) : g_L \leq u \leq g_U \textrm{ a.e.\ in } \Omega \right\}.
\end{equation*}
Several Schwarz methods for obstacle problems of the form~\eqref{obs_VI} were proposed in, e.g.,~\cite{Badea:2006,BTW:2003,Tai:2003,THX:2002}.
One can readily show that the variational inequality~\eqref{obs_VI} is equivalent to the following minimization problem:
\begin{equation}
\label{obs_opt}
\min_{u \in K} \frac{1}{2} \intO |\nabla u|^2 \,dx.
\end{equation}

Let $\Omega = [0,1]^2 \subset \mathbb{R}^2$.
We again use the discretization and domain decomposition settings introduced in Section~\ref{Subsec:sLap}.
That is, we consider a conforming discretization of~\eqref{obs_opt} on the continuous and piecewise linear finite element space $S_h (\Omega)$.
One more thing to do is to define an appropriate discretization $K_h$ of the set $K$; we simply set
\begin{equation*}
K_h = \left\{ u \in S_h (\Omega) : I_h g_L \leq u \leq I_h g_U \right\},
\end{equation*}
where $I_h$ is the nodal interpolation operator onto $S_h (\Omega)$.
Finally, the resulting discretization of~\eqref{obs_opt} is written as
\begin{equation}
\label{obs_FEM}
\min_{u \in K_h} \frac{1}{2} \intO |\nabla u|^2 \,dx.
\end{equation}
One may regard the constrained problem~\eqref{obs_FEM} as a nonsmooth unconstrained optimization problem.
More precisely, the discrete problem~\eqref{obs_FEM} is an instance of~\eqref{model} with
\begin{equation*}
V = S_h (\Omega), \quad F(u) = \frac{1}{2} \intO |\nabla u|^2 \,dx, \quad G(u) = \chi_{K_h}(u),
\end{equation*}
where $\chi_{K_h}$:~$V \rightarrow \overline{\mathbb{R}}$ is the characteristic function of $K_h$:
\begin{equation}
\label{chi}
\chi_{K_h} (u) = \begin{cases} 0 & \textrm{ if } u \in K_h, \\
\infty & \textrm{ if } u \not\in K_h. \end{cases}
\end{equation}
Clearly, $F$ satisfies the sharpness condition~\eqref{sharp} with $p = 2$.
Under the domain decomposition settings~\eqref{space_decomp} and~\eqref{space_decomp_2L}, the following convergence theorem for Algorithm~\ref{Alg:ASM} was presented in~\cite[Theorem~6.3]{Park:2020}.

\begin{proposition}
\label{Prop:obs}
In Algorithm~\ref{Alg:ASM} for~\eqref{obs_FEM}, there exists a positive constant $C$ independent of $h$, $H$, and $\delta$ such that
\begin{equation*}
E(u^{(n)}) - E(u^*) \leq \left( 1 - \frac{1}{2} \min \left\{ \tau, \frac{C}{1 + 1/\delta^2} \right\} \right)^n ( E (u^{(0)}) - E(u^*) )
\end{equation*}
for the one-level domain decomposition~\eqref{space_decomp} with $\tau \leq 1/4$ and
\begin{equation*}
E(u^{(n)}) - E(u^*) \leq \left( 1 - \frac{1}{2} \min \left\{ \tau, \frac{C}{\left( 1 + \log (H/h) \right) \left( 1 + (H/\delta)^2 \right)} \right\} \right)^n ( E (u^{(0)}) - E(u^*) )
\end{equation*}
for the two-level domain decomposition~\eqref{space_decomp_2L} with $\tau \leq 1/5$.
\end{proposition}

\begin{remark}
\label{Rem:obs}
There have been several notable existing results on the convergence rate of the unaccelerated additive Schwarz method for~\eqref{obs_FEM}~\cite{BTW:2003,Tai:2003,THX:2002}.
Proposition~\ref{Prop:obs} agrees with these existing results in the sense that the linear convergence rate is dependent on the stable decomposition property~(see, e.g.,~\cite[Eq.~(6)]{THX:2002} and~\cite[Assumption~4.1]{Park:2020}) of the space decomposition.
\end{remark}

For numerical experiments for the problem~\eqref{obs_FEM}, we set the obstacle functions $g_L$ and $g_U$ by
\begin{equation*}
g_L(x,y) = \begin{cases} 1 & \textrm{ if } \left(x - \frac{1}{2}\right)^2 + \left(y - \frac{1}{2} \right)^2 \leq \frac{1}{16^2}, \\
0 & \textrm{ otherwise,} \end{cases}
\end{equation*}
and
\begin{equation*}
g_U(x,y) = \begin{cases} 0 & \textrm{ if } \left(x - \frac{1}{4}\right)^2 + \left(y - \frac{1}{4} \right)^2 \leq \frac{1}{16^2}, \\
1 & \textrm{ otherwise,} \end{cases}
\end{equation*}
respectively.
The initial guess was chosen as $u^{(0)} = I_h g_L$ in order to make the condition $u^{(0)} \in \dom G$ in Algorithms~\ref{Alg:ASM} and~\ref{Alg:proposed} holds.
Local problems on $V_k$, $1\leq k \leq \N$, were solved by Algorithm~\ref{Alg:FISTA_ada} accompanied with backtracking~\cite{BT:2009} with the stop criterion~\eqref{stop_local}, while coarse problems on $V_0$ were solved by the nonlinear Gauss--Seidel method introduced in~\cite[section~5]{BTW:2003} with the stop criterion~\eqref{stop_coarse}.
The step size $\tau$ in Algorithms~\ref{Alg:ASM} and~\ref{Alg:proposed} was chosen in the same way as in Section~\ref{Subsec:sLap}.

\begin{figure}[]
\centering
\subfloat[][One-level decomposition~\eqref{space_decomp}]{ \includegraphics[width=0.47\linewidth]{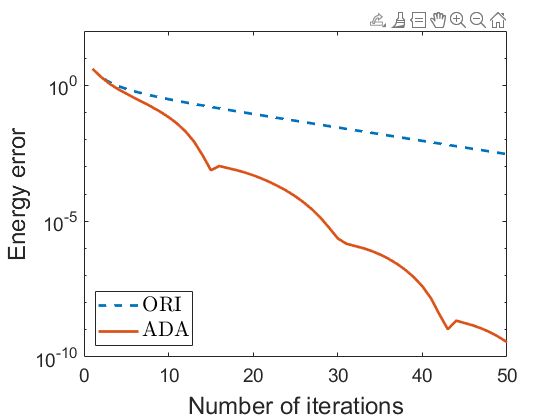} }
\gap
\subfloat[][Two-level decomposition~\eqref{space_decomp_2L}]{ \includegraphics[width=0.47\linewidth]{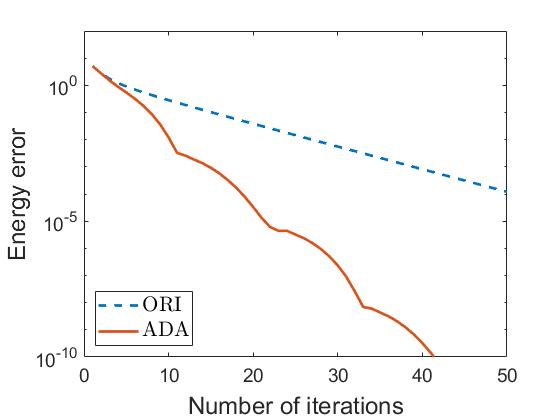} }
\caption{Decay of the energy error $E(u^{(n)}) - E(u^*)$ in additive Schwarz methods for the obstacle problem~\eqref{obs_FEM}~($h = 1/2^6$, $H=1/2^3$, $\delta=4h$).
ORI and ADA denote the unaccelerated~(Algorithm~\ref{Alg:ASM}) and accelerated~(Algorithm~\ref{Alg:proposed}) methods, respectively.}
\label{Fig:obs}
\end{figure}

In Figure~\ref{Fig:obs}, we present the energy error $E(u^{(n)})- E(u^*)$ of Algorithms~\ref{Alg:ASM} and~\ref{Alg:proposed} for~\eqref{obs_FEM} with respect to the number of iterations $n$,
where a reference solution $u^* \in V$ was obtained by $10^5$ iterations of Algorithm~\ref{Alg:FISTA_ada} applied to the full-dimension problem.
We note that Algorithm~\ref{Alg:ASM} for~\eqref{obs_FEM} is identical to~\cite[Algorithm~1]{THX:2002}.
It is observed that the convergence rate Algorithm~\ref{Alg:proposed} is faster than the one of Algorithm~\ref{Alg:ASM} for both one-level and two-level cases.
Therefore, we conclude that the proposed method is superior to the unaccelerated method when it is applied to the obstacle problem.

\begin{figure}[]
\centering
\subfloat[][$H/h = 2^3$, $\delta/h = 2$]{ \includegraphics[width=0.47\linewidth]{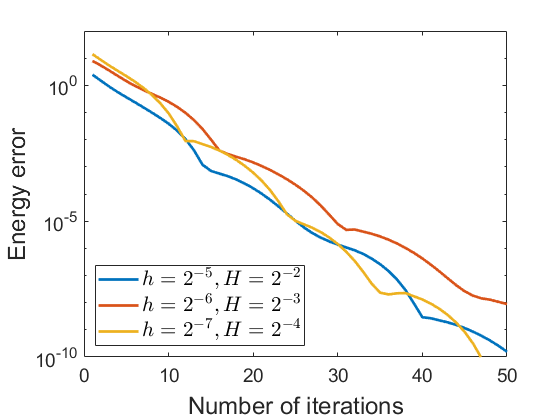} }
\gap
\subfloat[][$H/h = 2^4$, $\delta/h = 2$]{ \includegraphics[width=0.47\linewidth]{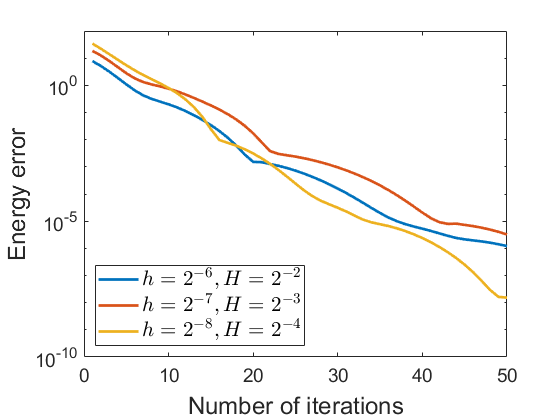} }
\caption{Decay of the energy error $E(\un) - E(u^*)$ in Algorithm~\ref{Alg:proposed} for the obstacle problem~\eqref{obs_FEM} when $H/h$ and $\delta/h$ are fixed.}
\label{Fig:obs_sca}
\end{figure}

\begin{table}
\centering
\begin{tabular}{c c c c} \hline
$H/h$ & $H$ & $h$ & \#iter \\
\Xhline{1pt}
\multirow{3}{*}{$2^2$} & $1/2^2$ & $1/2^4$ & 21 \\
& $1/2^3$ & $1/2^5$ & 35 \\
& $1/2^4$ & $1/2^6$ & 31  \\
\hline
\multirow{3}{*}{$2^3$} & $1/2^2$ & $1/2^5$ & 39 \\
& $1/2^3$ & $1/2^6$ & 50 \\
& $1/2^4$ & $1/2^7$ & 41 \\
\hline
\multirow{3}{*}{$2^4$} & $1/2^2$ & $1/2^6$ & 64 \\
& $1/2^3$ & $1/2^7$ & 72 \\
& $1/2^4$ & $1/2^8$ & 53 \\
\hline
\end{tabular}
\caption{Number of iterations of Algorithm~\ref{Alg:proposed} for the obstacle problem~\eqref{obs_FEM} to meet the condition~\eqref{stop_sca}.}
\label{Table:obs_sca}
\end{table}

Similarly to the case of the $s$-Laplacian problem, the scalability Algorithm~\ref{Alg:ASM} for the obstacle problem~\eqref{obs_FEM} is ensured by Proposition~\ref{Prop:obs}.
Therefore, the proposed method is also anticipated to enjoy the scalability; we do some numerical experiments for verification.
As shown in Figure~\ref{Fig:obs_sca}, the slopes of the energy graphs plotted in logarithmic scale in energy are almost indistinguishable for various values of $h$ and $H$ whenever the ratio $H/h$ is fixed.
Moreover, the number of iterations to meet the stop condition~\eqref{stop_sca} presented in Table~\ref{Table:obs_sca} does not increase even if both $H$ and $h$ decrease such that $H/h$ is kept constant.
These results verify that Algorithm~\ref{Alg:proposed} for~\eqref{obs_FEM} is numerically scalable.

\subsection{Dual total variation minimization}
Since it was numerically shown in~\cite[section~5.1]{OC:2015} that adaptive restarts can improve the performance of FISTA even for the nonsharp case, we expect that the proposed method may outperforms Algorithm~\ref{Alg:ASM} for the case without sharpness. 
As an example lacking the sharpness condition~\eqref{sharp}, we consider a minimization problem
\begin{equation}
\label{dualTV}
\min_{u \in H_0 (\div; \Omega)} \frac{1}{2} \intO ( \div u + f )^2 \,dx \quad
\textrm{ subject to } |u| \leq 1 \textrm{ a.e.\ in } \Omega,
\end{equation}
where $f \in L^2 (\Omega)$ and $|\cdot|$ denotes the pointwise Euclidean norm.
It is clear that~\eqref{dualTV} does not satisfy~\eqref{sharp} due to the divergence operator therein.
Problems of the form~\eqref{dualTV} usually appear in the field of mathematical imaging as dual formulations of total variation minimization problems~\cite{LPP:2019,LP:2019}.
Some overlapping Schwarz methods for~\eqref{dualTV} were studied in, e.g.,~\cite{CTWY:2015,HL:2015,Park:2021}.

Let $\Omega = [0,1] \subset \mathbb{R}^2$.
The domain $\Omega$ is decomposed into $\N = N \times N$ nonoverlapping subdomains $\{ \Omega_k \}_{k=1}^{\N}$, so that each subdomain has the sidelength $H = 1/N$.
Each subdomain $\Omega_k$ is partitioned into $H/h \times H/h$ uniform square elements.
Let $\Q_h$ be a subdivision of $\Omega$ consisting of those square elements.
We enlarge the subdomain $\Omega_k$ by adding several layers of square elements of width $\delta$ surrounding $\Omega_k$ to construct a region $\Omega_k'$, where $0< \delta < H/2$.
Then $\{ \Omega_k'\}$ forms an overlapping decomposition of $\Omega$.

We define $\bS_h (\Omega)$ as the lowest-order Raviart--Thomas finite element space on $\T_h$ with the homogeneous essential boundary condition.
In addition, let $K_h$ be a convex subset of $\bS_h (\Omega)$ defined by
\begin{equation*}
K_h = \left\{ u \in \bS_h (\Omega) : \frac{1}{|e|} \int_{e} |u \cdot \mathbf{n}_e| ds \leq 1, \gap e\textrm{: interior edges of } \Q_h \right\},
\end{equation*}
where $\mathbf{n}_e$ denotes the unit outer normal to $e$.
In~\cite{HHSVW:2019,LPP:2019}, the following discretization of~\eqref{dualTV} constructed by replacing the solution space and the constraint set by $\bS_h (\Omega)$ and $K_h$, respectively, was proposed:
\begin{equation}
\label{dualTV_FEM}
\min_{u \in K_h } \frac{1}{2} \intO ( \div u + f )^2 \,dx.
\end{equation}
Then the discrete problem~\eqref{dualTV_FEM} is of the form~\eqref{model} with
\begin{equation*}
V = \bS_h (\Omega), \quad F(u) = \frac{1}{2} \intO (\div u + f )^2 \,dx, \quad G(u) = \chi_{K_h}(u),
\end{equation*}
where $\chi_{K_h}$:~$V \rightarrow \overline{\mathbb{R}}$ is defined in the same manner as~\eqref{chi}.

In additive Schwarz methods for~\eqref{dualTV_FEM}, we set
\begin{equation*}
V_k = \bS_h (\Omega_k'), \quad 1\leq k \leq \N
\end{equation*}
an take $R_k^*$:~$V_k \rightarrow V$ as the natural extension operator so that~\eqref{space_decomp} holds, where $\bS_k (\Omega_k')$ is the lowest-order Raviart--Thomas finite element space on the $\Q_h$-elements in $\Omega_k'$ with the homogeneous essential boundary condition.
Then, according to~\cite[Theorem~6.5]{Park:2020}, we have the following convergence theorem.

\begin{proposition}
\label{Prop:dualTV}
In Algorithm~\ref{Alg:ASM} for~\eqref{dualTV_FEM}, there exists a positive constant $C$ independent of $h$, $H$, and $\delta$ such that
\begin{equation*}
E(u^{(n)}) - E(u^*) \leq C\frac{1 + 1/\delta^2}{n+1}
\end{equation*}
for the domain decomposition~\eqref{space_decomp} with $\tau \leq 1/4$.
\end{proposition}

Moreover, it was recently studied in~\cite{Park:2021} that the additive Schwarz method for~\eqref{dualTV_FEM} has a property called pseudo-linear convergence; it converges as fast as a linearly convergent algorithm until its energy error reaches a particular value which depends on $\delta$ only.
The following proposition summarizes the pseudo-linear convergence of the additive Schwarz method for~\eqref{dualTV_FEM}~\cite[Corollary~4.4]{Park:2021}.

\begin{proposition}
\label{Prop:dualTV_pseudo}
In Algorithm~\ref{Alg:ASM} for~\eqref{dualTV_FEM}, there exist two positive constants $\gamma <1$ and $\epsilon \leq C/\delta^2$ such that
\begin{equation*}
E(u^{(n)}) - E(u^*) \leq \gamma^n \left(E(u^{(0)}) - E(u^*) - \epsilon \right) + \epsilon
\end{equation*}
for the domain decomposition~\eqref{space_decomp} with $\tau \leq 1/4$, where $\gamma$ and $C$ are independent of $h$, $H$, and $\delta$.
\end{proposition}

For numerical experiments for the dual total variation minimization~\eqref{dualTV_FEM}, we set
\begin{equation*}
f(x,y) = \begin{cases} 1 & \textrm{ if } \left(x - \frac{1}{2}\right)^2 + \left(y - \frac{1}{2} \right)^2 \leq \frac{1}{4^2}, \\
0 & \textrm{ otherwise.} \end{cases}
\end{equation*}
The zero initial guess $u^{(0)} = 0$ were used.
We solved local problems on $V_k$, $1 \leq k \leq \N$, by Algorithm~\ref{Alg:FISTA_ada} with the step size $\tau = 1/8$~(see~\cite[Proposition~2.5]{LPP:2019}) and the stop criterion
\begin{equation*}
h^2 \| \div (w_k^{(n+1)} - w_k^{(n)}) \|_{\ell^2}^2 < 10^{-20}.
\end{equation*}
The step size for additive Schwarz methods was chosen as $\tau = 1/4$.
A reference solution $u^* \in V$ was obtained by Algorithm~\ref{Alg:FISTA_ada} with the step size $\tau=1/8$ for the full-dimension problem~\eqref{dualTV_FEM}.

\begin{figure}[]
\centering
 \includegraphics[width=0.47\linewidth]{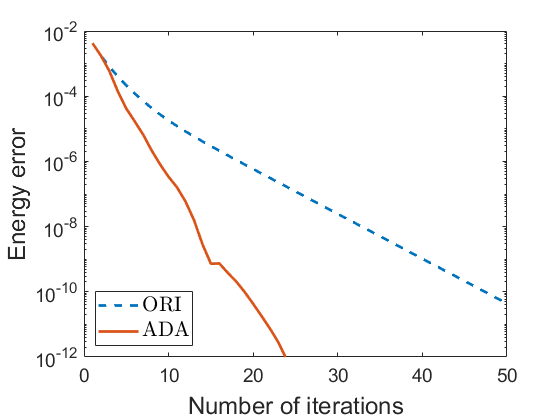}
\caption{Decay of the energy error $E(u^{(n)}) - E(u^*)$ in additive Schwarz methods for the dual total variation minimization~\eqref{dualTV_FEM}~($h = 1/2^6$, $H=1/2^3$, $\delta=4h$).
ORI and ADA denote the unaccelerated~(Algorithm~\ref{Alg:ASM}) and accelerated~(Algorithm~\ref{Alg:proposed}) methods, respectively.}
\label{Fig:dualTV}
\end{figure}

Figure~\ref{Fig:dualTV} plots $E(u^{(n)}) - E(u^*)$ of Algorithms~\ref{Alg:ASM} and~\ref{Alg:proposed} for~\eqref{dualTV_FEM} per iteration.
As in the cases of the $s$-Laplacian and obstacle problems, the proposed method outperforms Algorithm~\ref{Alg:ASM} in view of convergence rate.
Therefore, we conclude that the proposed method leads an improvement of the convergence rate even in the case of absence of the sharpness.

\begin{remark}
\label{Rem:dualTV_2L}
To the best of our knowledge, there have been no existing two-level Schwarz methods for dual total variation minimization problems, so that we do not include numerical results for the two-level case for~\eqref{dualTV_FEM} in this paper.
Nevertheless, we expect that if a suitable two-level domain decomposition method for~\eqref{dualTV_FEM} were developed in a near future,  the acceleration scheme presented would be applicable to that method to yield an improved convergence result.
\end{remark}

\begin{figure}[]
\centering
\subfloat[][$h = 1/2^7$, $\delta/h = 4$]{ \includegraphics[width=0.47\linewidth]{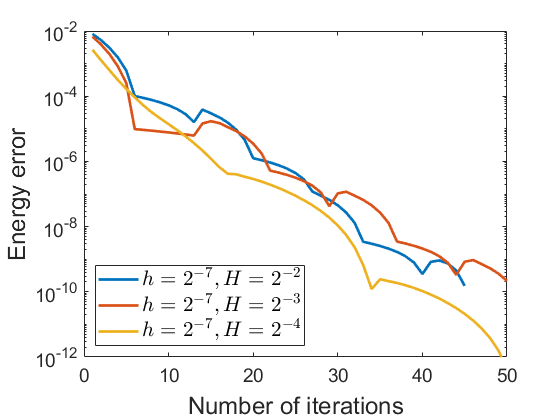} }
\gap
\subfloat[][$h = 1/2^8$, $\delta/h = 8$]{ \includegraphics[width=0.47\linewidth]{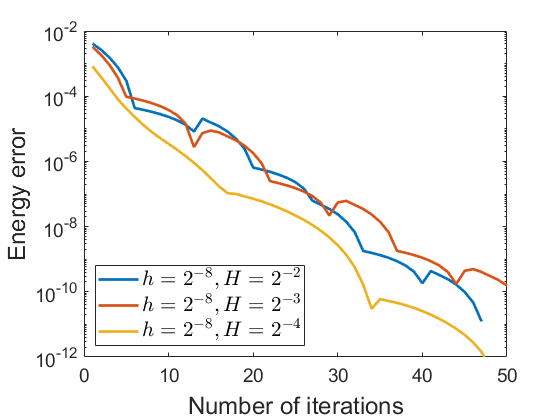} }
\caption{Decay of the energy error $E(\un) - E(u^*)$ in Algorithm~\ref{Alg:proposed} for the dual total variation minimization~\eqref{dualTV_FEM} when $\delta$ is fixed as $1/2^5$.}
\label{Fig:dualTV_pseudo}
\end{figure}

A remarkable property of Algorithm~\ref{Alg:ASM} for the dual total variation minimization~\eqref{dualTV_FEM} is pseudo-linear convergence; even though the unified theory of additive Schwarz methods for convex optimization presented in~\cite{Park:2020} gives only the sublinear convergence as stated in Proposition~\ref{Prop:dualTV}, it is verified by Proposition~\ref{Prop:dualTV_pseudo}  that Algorithm~\ref{Alg:ASM} for~\eqref{dualTV_FEM} performs as good as a linearly convergent algorithm if $\delta$ is large enough.
Therefore, we can expect that Algorithm~\ref{Alg:proposed} also enjoys pseudo-linear convergence since it was shown above by experiments that Algorithm~\ref{Alg:proposed} converges much faster than Algorithm~\ref{Alg:ASM}.
Indeed, Figure~\ref{Fig:dualTV_pseudo} shows how the energy error of Algorithm~\ref{Alg:proposed} for~\eqref{dualTV_FEM} decays when $\delta$ is fixed as $1/2^5$ while $H$ and $h$ vary.
It seems that all the cases share the approximately same linear convergence rate, similarly to Algorithm~\ref{Alg:ASM}; see~\cite{Park:2021} for the corresponding numerical results for Algorithm~\ref{Alg:ASM}.

\subsection{Comparison with the conjugate gradient method}
For completeness, we discuss the performance of the proposed method when it is applied to a linear elliptic problem.
It is well-known that the conjugate gradient method is optimal for symmetric and positive definite linear systems in the sense that it always find a solution in finite steps; see, e.g.,~\cite[Lemma~C.8]{TW:2005}.
Therefore, we may anticipate that the convergence rate of the proposed method lies between the ones of the Richardson method~(Algorithm~\ref{Alg:ASM}) and the conjugate gradient method.

In the following, we present numerical results of Algorithm~\ref{Alg:ASM}, Algorithm~\ref{Alg:proposed}, and the preconditioned conjugated gradient method with the additive Schwarz preconditioner~\eqref{preconditioner} for the problem~\eqref{sLap_FEM} with $s=2$ and $f = 1$.
Note that~\eqref{sLap_FEM} reduces the well-known Poisson equation when $s=2$.
Since local and coarse problems are linear, they can be solved directly, e.g., by the Cholesky factorization.
A reference solution $u^* \in V$ also can be obtained by a direct solver.

\begin{figure}[]
\centering
\subfloat[][One-level decomposition~\eqref{space_decomp}]{ \includegraphics[width=0.47\linewidth]{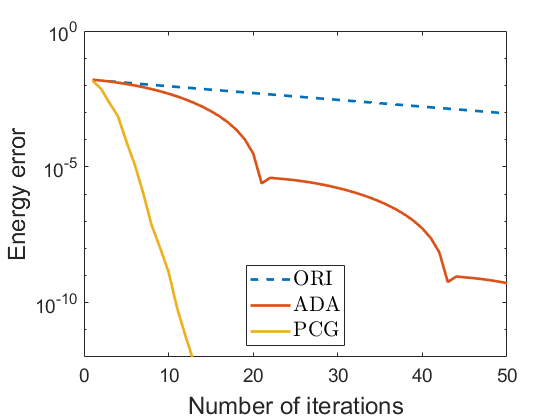} }
\gap
\subfloat[][Two-level decomposition~\eqref{space_decomp_2L}]{ \includegraphics[width=0.47\linewidth]{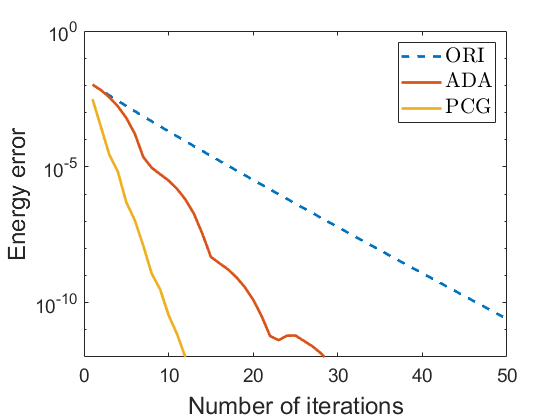} }
\caption{Decay of the energy error $E(u^{(n)}) - E(u^*)$ in additive Schwarz methods for the Poisson problem~($h = 1/2^6$, $H=1/2^3$, $\delta=4h$).
ORI and ADA denote the unaccelerated~(Algorithm~\ref{Alg:ASM}) and accelerated~(Algorithm~\ref{Alg:proposed}) methods, respectively.
PCG denotes the preconditioned conjugate gradient method with the additive Schwarz preconditioner~\eqref{preconditioner}.}
\label{Fig:Poisson}
\end{figure}

Figure~\ref{Fig:Poisson} depicts the convergence rates of Algorithm~\ref{Alg:ASM}, Algorithm~\ref{Alg:proposed}, and the corresponding preconditioned conjugate gradient method applied to the Poisson problem for both the one-level~\eqref{space_decomp} and two-level~\eqref{space_decomp_2L} cases.
As expected, the proposed method converges faster than Algorithm~\ref{Alg:ASM}, but slower than the preconditioned conjugate gradient method.
While there are several optimal solvers for linear problems such as the conjugate gradient method~(see~\cite{Saad:2003} for instance), designing optimal solvers for the general convex optimization of the form~\eqref{model} is still an important research topic that is actively investigated nowadays~\cite{DT:2014,KF:2016,RD:2020}.
In this perspective, the proposed method can be a very effective option for parallel computation of nonlinear and nonsmooth problems of the form~\eqref{model}.

\section{Conclusion}
\label{Sec:Conclusion}
In this paper, we proposed an accelerated additive Schwarz method that can be applied to very broad range of convex optimization problems of the form~\eqref{model}.
The proposed method showed superior convergence properties compared to the unaccelerated method when they are applied to various problems: $s$-Laplacian problem, obstacle problem, and dual total variation minimization.
Moreover, we verified that the proposed method inherits several desirable properties such as scalability and pseudo-linear convergence from the original method.
Since the proposed method does not require any a priori spectral information of a target problem, it has great potential to be applied to various convex optimization problems in science and engineering.

This paper leaves a few important topics for future research.
First, mathematical verification of faster convergence of the proposed method is required.
Unfortunately, to the best of our knowledge, there is no existing complete analysis even for Algorithm~\ref{Alg:FISTA_ada}, which is a main ingredient of the proposed method.
Meanwhile, we note that there is some recent research on acceleration of nonlinear gradient methods for the form~\eqref{GD}~\cite{Teboulle:2018}.
Another interesting topic is optimizing the acceleration scheme.
After a pioneering work~\cite{DT:2014}, there have been attempts to optimize acceleration schemes for gradient methods; see, e.g.,~\cite{KF:2016}.
We expect that it is possible to obtain a faster accelerated additive Schwarz method than the proposed method if we successfully apply such optimizing schemes to our case.

%
%
%

\bibliographystyle{spmpsci}
\bibliography{refs_ASM_restart}

\begin{thebibliography}{10}
\providecommand{\url}[1]{{#1}}
\providecommand{\urlprefix}{URL }
\expandafter\ifx\csname urlstyle\endcsname\relax
  \providecommand{\doi}[1]{DOI~\discretionary{}{}{}#1}\else
  \providecommand{\doi}{DOI~\discretionary{}{}{}\begingroup
  \urlstyle{rm}\Url}\fi

\bibitem{Badea:2006}
Badea, L.: Convergence rate of a {S}chwarz multilevel method for the
  constrained minimization of nonquadratic functionals.
\newblock SIAM J. Numer. Anal. \textbf{44}(2), 449--477 (2006)

\bibitem{BK:2012}
Badea, L., Krause, R.: One-and two-level {S}chwarz methods for variational
  inequalities of the second kind and their application to frictional contact.
\newblock Numer. Math. \textbf{120}(4), 573--599 (2012)

\bibitem{BTW:2003}
Badea, L., Tai, X.C., Wang, J.: Convergence rate analysis of a multiplicative
  {S}chwarz method for variational inequalities.
\newblock SIAM J. Numer. Anal. \textbf{41}(3), 1052--1073 (2003)

\bibitem{BT:2009}
Beck, A., Teboulle, M.: A fast iterative shrinkage-thresholding algorithm for
  linear inverse problems.
\newblock SIAM J. Imaging Sci. \textbf{2}(1), 183--202 (2009)

\bibitem{CD:2015}
Chambolle, A., Dossal, C.: On the convergence of the iterates of the ``{F}ast
  {I}terative {S}hrinkage/{T}hresholding {A}lgorithm''.
\newblock J. Optim. Theory Appl. \textbf{166}(3), 968--982 (2015)

\bibitem{CP:2016}
Chambolle, A., Pock, T.: An introduction to continuous optimization for
  imaging.
\newblock Acta Numer. \textbf{25}, 161--319 (2016)

\bibitem{CTWY:2015}
Chang, H., Tai, X.C., Wang, L.L., Yang, D.: Convergence rate of overlapping
  domain decomposition methods for the {R}udin--{O}sher--{F}atemi model based
  on a dual formulation.
\newblock SIAM J. Imaging Sci. \textbf{8}(1), 564--591 (2015)

\bibitem{Ciarlet:2002}
Ciarlet, P.G.: The Finite Element Method for Elliptic Problems.
\newblock SIAM, Philadelphia (2002)

\bibitem{CW:2005}
Combettes, P.L., Wajs, V.R.: Signal recovery by proximal forward-backward
  splitting.
\newblock Multiscale Model. Simul. \textbf{4}(4), 1168--1200 (2005)

\bibitem{DT:2014}
Drori, Y., Teboulle, M.: Performance of first-order methods for smooth convex
  minimization: a novel approach.
\newblock Math. Program. \textbf{145}(1-2), 451--482 (2014)

\bibitem{HHSVW:2019}
Herrmann, M., Herzog, R., Schmidt, S., Vidal-N{\'u}{\~n}ez, J., Wachsmuth, G.:
  Discrete total variation with finite elements and applications to imaging.
\newblock J. Math. Imaging Vision \textbf{61}(4), 411--431 (2019)

\bibitem{HL:2015}
Hinterm{\"u}ller, M., Langer, A.: Non-overlapping domain decomposition methods
  for dual total variation based image denoising.
\newblock J. Sci. Comput. \textbf{62}(2), 456--481 (2015)

\bibitem{KF:2016}
Kim, D., Fessler, J.A.: Optimized first-order methods for smooth convex
  minimization.
\newblock Math. Program. \textbf{159}(1-2), 81--107 (2016)

\bibitem{LPP:2019}
Lee, C.O., Park, E.H., Park, J.: A finite element approach for the dual
  {R}udin--{O}sher--{F}atemi model and its nonoverlapping domain decomposition
  methods.
\newblock SIAM J. Sci. Comput. \textbf{41}(2), B205--B228 (2019)

\bibitem{LP:2019b}
Lee, C.O., Park, J.: Fast nonoverlapping block {J}acobi method for the dual
  {R}udin--{O}sher--{F}atemi model.
\newblock SIAM J. Imaging Sci. \textbf{12}(4), 2009--2034 (2019)

\bibitem{LP:2019}
Lee, C.O., Park, J.: A finite element nonoverlapping domain decomposition
  method with {L}agrange multipliers for the dual total variation
  minimizations.
\newblock J. Sci. Comput. \textbf{81}(3), 2331--2355 (2019)

\bibitem{Nesterov:2013}
Nesterov, Y.: Gradient methods for minimizing composite functions.
\newblock Math. Program. \textbf{140}(1), 125--161 (2013)

\bibitem{Nesterov:2018}
Nesterov, Y.: Lectures on Convex Optimization.
\newblock Springer, Cham (2018)

\bibitem{Nesterov:1983}
Nesterov, Y.E.: A method for solving the convex programming problem with
  convergence rate ${O}(1/k^2)$.
\newblock Dokl. Akad. Nauk SSSR \textbf{269}, 543--547 (1983)

\bibitem{OC:2015}
O’Donoghue, B., Candes, E.: Adaptive restart for accelerated gradient
  schemes.
\newblock Found. Comput. Math. \textbf{15}(3), 715--732 (2015)

\bibitem{Park:2020}
Park, J.: Additive {S}chwarz methods for convex optimization as gradient
  methods.
\newblock SIAM J. Numer. Anal. \textbf{58}(3), 1495--1530 (2020)

\bibitem{Park:2021}
Park, J.: Pseudo-linear convergence of an additive {S}chwarz method for dual
  total variation minimization.
\newblock Electron. Trans. Numer. Anal. \textbf{54}, 176--197 (2021)

\bibitem{RD:2020}
Roulet, V., d'Aspremont, A.: Sharpness, restart, and acceleration.
\newblock SIAM J. Optim. \textbf{30}(1), 262--289 (2020)

\bibitem{Saad:2003}
Saad, Y.: Iterative Methods for Sparse Linear Systems.
\newblock SIAM, Philadelphia (2003)

\bibitem{Tai:2003}
Tai, X.C.: Rate of convergence for some constraint decomposition methods for
  nonlinear variational inequalities.
\newblock Numer. Math. \textbf{93}(4), 755--786 (2003)

\bibitem{THX:2002}
Tai, X.C., Heimsund, B., Xu, J.: Rate of convergence for parallel subspace
  correction methods for nonlinear variational inequalities.
\newblock In: Domain decomposition methods in science and engineering ({L}yon,
  2000), Theory Eng. Appl. Comput. Methods, pp. 127--138. Internat. Center
  Numer. Methods Eng. (CIMNE), Barcelona (2002)

\bibitem{TX:2002}
Tai, X.C., Xu, J.: Global and uniform convergence of subspace correction
  methods for some convex optimization problems.
\newblock Math. Comp. \textbf{71}(237), 105--124 (2002)

\bibitem{Teboulle:2018}
Teboulle, M.: A simplified view of first order methods for optimization.
\newblock Math. Program. \textbf{170}(1), 67--96 (2018)

\bibitem{TW:2005}
Toselli, A., Widlund, O.: Domain Decomposition Methods---Algorithms and Theory.
\newblock Springer, Berlin (2005)

\bibitem{Xu:1992}
Xu, J.: Iterative methods by space decomposition and subspace correction.
\newblock SIAM Rev. \textbf{34}(4), 581--613 (1992)

\end{thebibliography}

\end{document}